\documentclass{ifacconf}

\usepackage{amsmath,amsfonts,amssymb}
\usepackage{mathtools}
\usepackage{enumerate}
\usepackage{color}
\usepackage{tikz-cd} 
\usepackage{circuitikz}
\usepackage{pgfplots}
\usepackage{pgfplotstable}
\usepgfplotslibrary{patchplots}
\usepgfplotslibrary{groupplots}

\pgfplotsset{compat=1.5}

\usepackage{graphicx}      
\usepackage{natbib}        

\DeclareMathOperator{\pos}{pos}
\DeclareMathOperator{\interior}{int}
\DeclareMathOperator{\boundary}{bdr}

\DeclareMathOperator{\co}{co}

\DeclareMathOperator{\proj}{Proj}

\newcommand\RE{\mathbb{R}}
\newcommand\Fi{\mathcal{F}}
\newcommand\Ge{\mathcal{G}}

\newcommand{\bscdot}{\boldsymbol{\cdot}}

\begin{document}
	
\begin{frontmatter}
\title{
A notion of equivalence for linear complementarity problems with application to the design of non-smooth bifurcations
}

\author[Fst]{Fernando Casta\~nos} 
\author[Scd]{Felix A. Miranda-Villatoro} 
\author[Thrd]{Alessio Franci}

\address[Fst]{Automatic Control Department, Cinvestav-IPN. Av. Instituto Polit\'ecnico Nacional 2508, 07360, CDMX, Mexico. email: fcastanos@ctrl.cinvestav.mx}
\address[Scd]{Department of Engineering, University of Cambridge. Trumpington Street, CB2 1PZ, Cambridge, UK. email: fam48@cam.ac.uk}
\address[Thrd]{Department of Mathematics, Universidad Nacional Aut\'onoma de M\'exico, Circuito Exterior S/N, C.U., 04510, CDMX, Mexico. email: afranci@ciencias.unam.mx}

\begin{abstract}                
Many systems of interest to control engineering can be modeled by linear complementarity problems. We introduce a new notion of equivalence between linear complementarity problems that sets the basis to translate the powerful tools of smooth bifurcation theory to this class of models. Leveraging this notion of equivalence, we introduce new tools to analyze, classify, and design non-smooth bifurcations in linear complementarity problems and their interconnection.
\end{abstract}


\end{frontmatter}
\section{Introduction}
\label{section:introduction}

Bifurcation theory is one of the most 
successful tools for the analysis of nonlinear dynamical systems that depend on a 
control parameter. The theory is firmly grounded on the classical implicit function 
theorem~\citep{dontchev2014,golubitsky1985}, and therefore, it requires smoothness of the maps under study.
However, from a practical viewpoint, it is common to approximate complicated 
nonlinear maps by simpler models. In such situations, 
the resulting approximation may be non-smooth. 

Linear complementarity problems are non-smooth problems that arise in 
fields of science such as economics \citep{nagurney1993}, electronics \citep{acary2011}, 
mechanics \citep{brogliato1999}, 
mathematical programming \citep{murty1988}, general systems theory~\citep{schaft1998}, etc. 
They serve as a departing point in the analysis of problems with unilateral constraints,
and also arise as piecewise linear approximations of nonlinear models \citep{leenaerts1998}. 

Recently, there have been some attempts to extend bifurcation theory towards the 
non-smooth setting, see e.g. \cite{diBernardo2008, leine2004, simpson2010}.
However, the emphasis has been directed towards analysis of 
discontinuous systems, and very little is known on bifurcations in complementarity 
systems.

The purpose of this paper is to provide a methodology for the realization of equilibrium bifurcations 
in linear complementarity problems. The proposed framework mimics, up to certain extent, the smooth program proposed by \cite{arnold} and relies on tools from 
non-smooth analysis and linear algebra. To achieve this, the concept of topological equivalence in complementarity systems is introduced. We focus on static models that arise as the steady-state equations of piecewise linear dynamical systems. Thanks to the piecewise linear structure 
of the problem, the introduced equivalence is always global, which constitutes a major difference with respect to smooth bifurcation theories. This fundamental concept allows us to 
provide a complete classification of planar complementarity problems.

The paper is organized as follows. Section \ref{section:preliminaries} describes the linear complementarity problem and related concepts. Section \ref{section:main:results} constitutes
the main body of the paper and addresses the problem of topological equivalence between LCP's. Afterwards, an interconnection approach for the realization of bifurcations is presented, together with an example applied to the non-smooth pleat and the pitchfork singularity. 
Finally, the paper ends with some conclusions and future research directions in Section \ref{section:discussion:future}.


\section{Preliminaries}
\label{section:preliminaries}


\subsection{Linear Complementarity Problems}

The linear complementarity problem (LCP) is defined as follows.
\begin{defn}
  \label{defn:lcp}
 Given a vector $q \in \RE^n$ and a matrix $M \in \RE^{n\times n}$, the LCP $(M,q)$ consists in finding vectors $z,w \in \RE^n$ such that
\begin{equation}
  \label{eq:lcp}
  \begin{cases}
    w = M z + q 
    \\
    \RE_{+}^{n} \ni w \perp z \in \RE_{+}^{n}
  \end{cases}
\end{equation}
where the second relation, called the \emph{complementary condition}, is the short form of the following three conditions: $w \in \RE_{+}^{n}$, $z \in \RE_{+}^{n}$, and $w^{\top} z = 0$.
\end{defn}

In what follows, we introduce some concepts that will be useful for studying the geometric 
structure of LCPs. 
Given $M$ and an index set $\alpha \subseteq \left\{ 1,\dots,n \right\}$, we define the
\emph{complementary matrix} $C_M(\alpha)$ as
\begin{displaymath}
 C_M(\alpha)_{\bscdot j} = 
  \begin{cases}
   -M_{\bscdot j} &\quad \text{if $j \in \alpha$} \\
    I_{\bscdot j} &\quad \text{if $j \not\in \alpha$}
  \end{cases} \;,
\end{displaymath}
where the subscript $\bscdot j$ denotes the $j$-th column. Now define the piecewise-linear function
\begin{equation} \label{eq:fM}
 f_M(x) = C_{-M}(\alpha) x \;, \quad x \in \pos C_I(\alpha) \;,
\end{equation}
where $\pos C_I(\alpha)$ is the cone generated by the columns of $C_I(\alpha)$. Note that the 
cones $\pos C_I(\alpha)$ are simply the $2^n$ orthants in $\RE^n$ indexed by $\alpha \subseteq \{1, \dots, n\}$, and that 
\begin{displaymath}
 f_M(\pos C_I(\alpha)) = \pos C_M(\alpha) \;.
\end{displaymath}

\begin{prop}[\cite{cottle}]
  \label{prop:lcp:pwl:equivalence}
  Let  $z \in \RE^{n}$ be a solution of the LCP $(M,q)$, then 
  $x = w - z \in \RE^n$ is a solution of
 \begin{equation} \label{eq:fMq}
  f_M(x) = q \;.
 \end{equation}
 Conversely, let $x \in \RE^{n}$ be a solution of \eqref{eq:fMq}, then 
 $z = \proj_{\RE_{+}^{n}}(-x) \in \RE_{+}^{n}$ is a solution of the LCP $(M,q)$. 
\end{prop}

Henceforth, we treat the LCP $(M, q)$ and \eqref{eq:fMq} as identical problems, in the sense that we only need to know the solution of one of them in order to know the solution of the other.


The solutions of the LCP $(M,q)$ depend on the geometry of the \emph{complementary cones} $\pos C_M(\alpha)$.
More precisely, there exists at least one solution $x$ of~\eqref{eq:fMq} 
for every $\alpha$ such that $q \in \pos C_M(\alpha)$. If $C_{M}(\alpha)$ is 
nonsingular the solution is unique, whereas there exists a continuum of solutions if $C_{M}(\alpha)$ is singular. Thus, for a given $q$, 
there can be no solutions, there can be one solution, multiple isolated solutions, or a continuum of solutions, 
depending on how many complementary cones $q$ belongs to and their properties.

\subsection{Bifurcations in LCPs}

In practical applications, the vector $q$ depends on a control, or \emph{bifurcation} parameter $\lambda\in\RE$. The bifurcation parameter can be an applied voltage or current in electronic circuits, a force or a torque in a mechanical system, or the amount of capital injection in an economic system. The goal of bifurcation theory is to understand how the number of 
solutions changes as the bifurcation parameter is varied. In LCPs we let
$
q=\bar q(\lambda),
$
where $\bar q:\RE\to\RE^n$ is at least continuous, although more regularity constraints can be imposed as needed. The mapping $\bar q$ defines a continuous curve, or \emph{path} in $\RE^n$. As $\lambda$ lets $q$ move along this path, the number of solution to the LCPs might change. Points where the number of solutions change are called \emph{bifurcation points}. 

\begin{exmp}
  \label{exmp:lcp:1}
Let us illustrate this idea in the simple case where the path is a line segment joining two distinct points $q_{i} \in \RE^{2}$, $i \in \{0, 1\}$, that is,
\[
  \bar q(\lambda)= (1 - \lambda) q_0 + \lambda q_1, \quad  \lambda \in [0, 1].
\]
In addition, let us set the matrix $M$ as 
\begin{equation}
  \label{eq:M:lcp:1}
  M = 
  \begin{bmatrix}
    1 & 2
    \\
    2 & 1
  \end{bmatrix}
\end{equation}
and proceed to analyze the two cases shown in Fig.~\ref{fig:lcp1}.

\begin{figure}
  \centering
  \includegraphics[width=0.45\columnwidth]{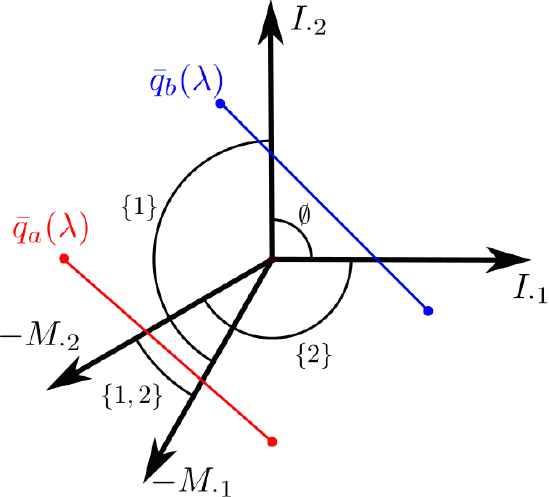}
  \caption{Cone configuration for matrix $M$ in \eqref{eq:M:lcp:1}. The thick black lines depict the generators of the complementary cones $\pos C_{M}(\alpha)$, $\alpha \subseteq \{1, 2 \}$, whereas the arcs denote the complementary cones.}
  \label{fig:lcp1}
\end{figure}

\emph{Case a)}
We take the path $\bar{q}_{a}(\lambda)$ given by
    \begin{equation}
      \label{eq:lcp1:path:a}
  \bar{q}_{a}(\lambda) = (1 - \lambda) 
  \begin{bmatrix}
    -4 \\ 0
  \end{bmatrix} + \lambda
  \begin{bmatrix}
    0 \\ -4
  \end{bmatrix} \;, \quad \lambda \in [0, 1] \;.
\end{equation}
According to Proposition \ref{prop:lcp:pwl:equivalence}, solving the LCP 
$(M, \bar{q}_{a}(\lambda))$ is equivalent to finding $x \in \RE^{2}$ satisfying 
\begin{equation}
  C_{-M}(\alpha) x = \bar{q}_{a}(\lambda), \quad \text{ s.t. } x \in \pos C_{I}(\alpha) \;,
  \label{eq:lcp1:pwl}
\end{equation}
for $\alpha \subseteq \{1, 2\}$.
Noting that $C_{M}(\alpha) = C_{-M}(\alpha) C_{I}(\alpha)$ for any $\alpha \subseteq \{1, 2\}$, it follows that the solutions to \eqref{eq:lcp1:pwl} are given by 
\begin{equation}
  \label{eq:lcp:sol:set}
  {\bigcup}_{\alpha \subseteq \{1, 2 \}} S_{\alpha}
\end{equation}
where 
\begin{multline}
  S_\alpha = \big\{ (x,\lambda) \in \RE^{2}\times[0,1] \mid \exists \; p_{\lambda}(\alpha) \in \RE_{+}^{2} \; : \\
   x = C_{I}(\alpha) p_{\lambda}(\alpha) \text{ and } \bar{q}_{a}(\lambda) = C_{M}(\alpha)p_{\lambda}(\alpha)
  \big\}
  \label{eq:lcp1:solution}
\end{multline}
%
Roughly speaking, in order to solve the parametrized LCP $(M, \bar{q}(\lambda))$ 
we need to find $p_{\lambda}(\alpha)$ (the representation of $\bar{q}(\lambda)$ 
in terms of the generators of the $\alpha$-th complementary cone).
Computing these explicitly and taking $\alpha = \emptyset \subset \{1, 2\}$ we get
\begin{displaymath}
  p_{\lambda}(\emptyset) = C_{M}(\emptyset)\bar{q}_{a}(\lambda) = 
  \begin{bmatrix}
    4 \lambda - 4 \\ -4 \lambda
  \end{bmatrix} \;,
\end{displaymath}
and it follows that $p_{\lambda}(\emptyset) \notin \RE_{+}^{2}$ for any $\lambda \in \RE$.
Therefore, $S_{\emptyset} = \emptyset$.
Now, for $\alpha = \{1\} \subset \{1, 2\}$ we get that
\begin{displaymath}
  p_{\lambda}(\{1\}) = C_{M}(\{1\}) \bar{q}_{a}(\lambda) = 
  \begin{bmatrix}
    4 - 4\lambda
    \\
    8 - 12 \lambda
  \end{bmatrix}.
\end{displaymath}
It follows that $p_{\lambda}(\{1\}) \in \RE_{+}^{2}$ for $\lambda \in (-\infty, 2/3]$.
Hence, 
\begin{displaymath}
  S_{\{1\}} = \left\{ (x, \lambda) \in \RE^{2} \times  [0, 2/3] \mid x  = 
     \begin{bmatrix}
      4 \lambda - 4 \\ 8 - 12 \lambda
  \end{bmatrix}\right\}
\end{displaymath}
Similarly, for $\alpha = \{2\}$ and $\alpha = \{1, 2\}$ we have, respectively,
\begin{align*}
  S_{\{2\}} &= \left\{ (x, \lambda) \in \RE^{2} \times [1/3, 1] \mid x = \begin{bmatrix}
      -4 + 12 \lambda \\ -4 \lambda
  \end{bmatrix}\right\} \;,
  \\
  S_{\{1, 2\}} &= \left\{ (x, \lambda) \in \RE^{2} \times [1/3, 2/3] 
    \mid x = \begin{bmatrix}
      \frac{4}{3} - 4 \lambda \\ 4 \lambda - \frac{8}{3}
  \end{bmatrix} \right\} \;. 
\end{align*}
By putting all the pieces together, one gets the bifurcation diagrams 
(in the $x$-coordinate) shown on the left-hand side of Fig.~\ref{fig:lcp1:x1}.

\emph{Case b)} We take the path
  \begin{equation}
    \label{eq:lcp1:path:b}
  \bar{q}_{b}(\lambda) = (1 - \lambda) 
  \begin{bmatrix}
    -1 \\ 3
  \end{bmatrix} + \lambda
  \begin{bmatrix}
    3 \\ -1
  \end{bmatrix}.
\end{equation}
As in the previous case, we need to solve a family of constrained linear problems. Simple computations lead us to
\begin{align*}
  S_{\emptyset} &= \left\{ (x, \lambda) \in \RE^{2} \times  [1/4, 3/4] \mid x = \begin{bmatrix}
    4 \lambda - 1 \\ 3 - 4 \lambda
\end{bmatrix} \right\}, \\
S_{\{1\}} &= \left\{ (x, \lambda) \in \RE^{2} \times [0, 1/4] \mid x = \begin{bmatrix}
    1 - 4 \lambda \\ 5 - 12 \lambda
\end{bmatrix} \right\}, 
\\
S_{\{2\}} &= \left\{ (x, \lambda) \in \RE^{2} \times [3/4, 1] \mid x = \begin{bmatrix}
  12 \lambda - 7 \\ 3 - 4 \lambda
\end{bmatrix} \right\}, 
\\
S_{\{1, 2\}} &= \emptyset \, .
\end{align*}

The right-hand side of Fig.~\ref{fig:lcp1:x1} depicts the  
solution set of LCP $(M, \bar{q}_{b}(\lambda))$ (in the $x$-variable).

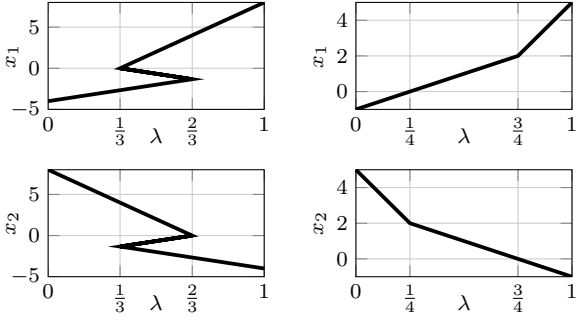
\begin{figure}
  \centering
  \small
  \begin{tikzpicture}
  \begin{groupplot}[group style={group size=1 by 2, vertical sep=8mm}, 
    height=3cm, width=0.5\columnwidth]
\nextgroupplot[
    xtick = {0.0, 0.333, 0.667, 1.0},
    xticklabels = {$0$, $\frac{1}{3}$, $\frac{2}{3}$, $1$},
x grid style={white!80.0!black},
y grid style={white!80.0!black},
xmajorgrids,
ymajorgrids,
x label style={at={(axis description cs:0.5,-0.1)},anchor=north},
y label style={at={(axis description cs:-0.1,0.5)},anchor=south},
xlabel={$\lambda$},
ylabel={$x_{1}$},
xmin=-0.0, xmax=1.0,
ymin=-5, ymax=8,
]
\addplot [line width=0.5mm, black, forget plot] table [x=l, y=x1, col sep=comma]{lcp1_hysteresis.csv};

\nextgroupplot[
    xtick = {0.0, 0.333, 0.667, 1.0},
    xticklabels = {$0$, $\frac{1}{3}$, $\frac{2}{3}$, $1$},
x grid style={white!80.0!black},
y grid style={white!80.0!black},
xmajorgrids,
ymajorgrids,
x label style={at={(axis description cs:0.5,-0.1)},anchor=north},
y label style={at={(axis description cs:-0.1,0.5)},anchor=south},
xlabel={$\lambda$},
ylabel={$x_{2}$},
xmin=-0.0, xmax=1.0,
ymin=-5, ymax=8,
]
\addplot [line width=0.5mm, black, forget plot] table [x=l, y=x2, col sep=comma]{lcp1_hysteresis.csv};

\end{groupplot}

\end{tikzpicture} \;
  \begin{tikzpicture}
\begin{groupplot}[group style={group size=1 by 2, vertical sep=8mm}, 
    height=3cm, width=0.5\columnwidth]
\nextgroupplot[
    xtick = {0.0, 0.25, 0.75, 1.0},
    xticklabels = {$0$, $\frac{1}{4}$, $\frac{3}{4}$, $1$},
x grid style={white!80.0!black},
y grid style={white!80.0!black},
xmajorgrids,
ymajorgrids,
x label style={at={(axis description cs:0.5,-0.1)},anchor=north},
y label style={at={(axis description cs:-0.1,0.5)},anchor=south},
xlabel={$\lambda$},
ylabel={$x_{1}$},
xmin=-0.0, xmax=1.0,
ymin=-1, ymax=5,
]
\addplot [line width=0.5mm, black, forget plot] table [x=l, y=x1, col sep=comma]{lcp1_regular.csv};

\nextgroupplot[
    xtick = {0.0, 0.25, 0.75, 1.0},
    xticklabels = {$0$, $\frac{1}{4}$, $\frac{3}{4}$, $1$},
x grid style={white!80.0!black},
y grid style={white!80.0!black},
xmajorgrids,
ymajorgrids,
x label style={at={(axis description cs:0.5,-0.1)},anchor=north},
y label style={at={(axis description cs:-0.1,0.5)},anchor=south},
xlabel={$\lambda$},
ylabel={$x_{2}$},
xmin=-0.0, xmax=1.0,
ymin=-1, ymax=5,
]
\addplot [line width=0.5mm, black, forget plot] table [x=l, y=x2, col sep=comma]{lcp1_regular.csv};
\end{groupplot}

\end{tikzpicture}
  \caption{Solutions to problem \eqref{eq:lcp1:pwl}
    for $M$ given by \eqref{eq:M:lcp:1} and paths
  $\bar{q}_{a}$ as in \eqref{eq:lcp1:path:a} (left), and $\bar{q}_{b}$ as in 
\eqref{eq:lcp1:path:b} (right), both with $\lambda \in [0, 1]$.}
  \label{fig:lcp1:x1}
\end{figure}
\end{exmp}

It is clear that, as long as the path $\bar{q}(\lambda)$ lies in the interior of the same cone, or set of cones, the number of solutions cannot change. 
Exiting and/or entering a cone, that is, crossing a cone face is thus a necessary condition for a bifurcation to occur. It is not sufficient though. 
For instance, in Example \ref{exmp:lcp:1} Case b) above, the path $\bar{q}_{b}(\lambda)$ crosses through different cones at the points $\lambda \in \{ \frac{1}{4}, \frac{3}{4} \}$. However, there is no change in the number of solutions, see  Fig.~\ref{fig:lcp1:x1}, right. This last observation poses the following question: How can we characterize the face at which bifurcations occur?


The non-smooth Implicit Function Theorem (see Corollary at page 256 of \cite{clarke1990}) provides an answer to this question. Let $\Omega_f$ be the set of measure zero where the Jacobian $Df(x)$ of a Lipschitz continuous function $f: \RE^n \to \RE^n$ does not exist.

\begin{defn}[Clarke generalized Jacobian]
 The generalized Jacobian of $f$ at $x$ is the set
 \[
  \partial f(x) = \co \left\{ \lim_{i \to \infty} D f(x_i) \mid x_i \to x, \; x_i \not\in S, \; x_i \not\in \Omega_f \right\} \;,
 \]
 where $S$ is any set of measure zero and $\co$ denotes convex closure.
\end{defn}

\begin{defn}
 $\partial f(x)$ is said to be \emph{maximal rank} if every $M$ in $\partial f(x)$ is non-singular.
\end{defn}

For a function $F:\RE^n\times\RE^m\to\RE^n$, $F:(x,y)\to F(x,y)$, the generalized Jacobian with respect to the first argument, denoted by $\partial_x F(x,y)$, is the set of all $n\times n$ matrices $M$ such that $\begin{bmatrix} M & N\end{bmatrix}$ belongs to $\partial F(x,y)$ for some $n\times m$ matrix $N$.

\begin{thm}
Suppose that $F(x_0,y_0)=0$ and its generalized Jacobian $\partial_x F(x_0,y_0)$ is maximal rank. Then there exist a neighborhood $U$ of $y_0$ and a Lipschitz function $\bar x:U\to\RE^n$ such that $F(\bar x(y),y)\equiv 0$ for all $y\in U$.
\end{thm}

By specializing this theorem to~\eqref{eq:fMq} with $F(x,q) = f_M(x)-q$, it follows that a solution $(x_0,q_0)$ to an LCP can be a bifurcation point only if $\partial f_M(x^*)$ is not maximal rank, that is, if there exists a singular matrix $M_0$ belonging to the set $\partial f_M(x^*)$. This motivates the following definition.

\begin{defn}
  \label{defn:nonsmooth:singularity}
A solution point $(x_0,q_0)$ of~\eqref{eq:fMq} such that $\partial f_M(x_0)$ is not maximal rank is called a non-smooth singularity.
\end{defn}

Observe that 
$ \partial f_M(x)=\co\left\{ C_{-M}(\alpha)\mid x\in\pos C_I(\alpha) \right\} $.
Thus, $\partial f_M(x)$ is a singleton if $x$ belongs to the interior of an orthant 
or the convex closure of a (finite) set of matrices if $x$ belongs to the face between two or more orthants.


%

The following proposition helps in finding non-smooth singular points.

\begin{prop}\label{prop:singular:sufficient}
Let $x_0$ be a solution of the LCP $(M,q_0)$. If there exists $M_+\in \partial f_M(x_0)$ such that $\det(M_+)>0$ and 
$M_{-}\in \partial f(x_0)$ such that $\det(M_{-})<0$, then $\partial f_M(x_0)$ is not maximal rank.
\end{prop}

\begin{pf}
The determinant function $\det:\RE^{n\times n}\to\RE$ is continuous and the set $\partial f_M(x_0)$ is connected (since it is convex). It follows that, because $\det(M)$ takes both positive and negative values in $\partial f_M(x_0)$, it must also vanish in some subset of $\partial f_M(x_0)$.
 \qed
\end{pf}


As an application of Proposition \ref{prop:singular:sufficient}, 
let us consider Example \ref{exmp:lcp:1} above. Note that $\partial f_{M}$ can be set-valued 
only for $\{ x \in \RE^{n} \mid f_{M}(x) \in \boundary \pos C_{M}(\alpha)$,
$\alpha \subseteq \{1, \dots, n \}\}$.
Therefore, with $M$ as in \eqref{eq:M:lcp:1}, the generalized Jacobian (at the coordinate axes) is 
\begin{equation}
 \partial f_{M}(x) = 
\begin{cases}
\left\{ \begin{bmatrix}
  1 & 0 \\ 2 - 2 \mu & 1
\end{bmatrix},  \mu \in [0, 1] \right\},
\\
\qquad \qquad \qquad x \in \pos C_{I}(\emptyset) \cap \pos C_{I}(\{1\})
\\
\left\{ \begin{bmatrix}
  1 & 2 - 2\mu \\ 0 & 1
\end{bmatrix},  \mu \in [0, 1] \right\}, 
\\
\qquad \qquad \qquad x \in \pos C_{I}(\emptyset) \cap \pos C_{I}(\{2\})
\\
\left\{
\begin{bmatrix}
  1 & 2 \mu \\ 2 & 1
\end{bmatrix},  \mu \in [0, 1] \right\},
\\
\qquad \qquad x \in \pos C_{I}(\{1,2\}) \cap \pos C_{I}(\{1\})
\\
\left\{
\begin{bmatrix}
  1 & 2 \\ 2 \mu & 1
\end{bmatrix}, \mu \in [0, 1] \right\},
\\
\qquad \qquad x \in \pos C_{I}(\{1,2\}) \cap \pos C_{I}(\{2\})
\end{cases}
  \label{eq:lcp1:jacobian}
\end{equation}
whereas $\partial f_{M}$ is single-valued and nonsingular for all of the others points
$x \in \RE^{2}$. 
It follows from Proposition \ref{prop:singular:sufficient} and \eqref{eq:lcp1:jacobian}
that solutions of $f_{M}(x) - q = 0$ satisfying 
$x_{1} = 0$ or $x_{2} = 0$ are non-smooth singular points,
see the expression for $S_{\{1, 2\}}$ and the left-hand side of Fig.~\ref{fig:lcp1:x1}.
In contrast, it follows directly from
Definition \ref{defn:nonsmooth:singularity} that 
all solutions 
of Case $b)$ in Example \ref{exmp:lcp:1} are regular, see the right-hand side of
Fig.~\ref{fig:lcp1:x1}.
It is worth to remark that to have a singularity it is not necessary that $\det (C_{-M}(\alpha))=0$ for some $\alpha$.

%
%
%
%
%
%
When $\det(C_{-M}(\alpha))=0$ for some $\alpha$ such that $q_0\in\pos C_{-M}(\alpha)$,
another source of singularities appears. 
In this case, the cone $\pos C_{-M}(\alpha)$ is \emph{degenerate}, in the sense that its $n$-dimensional interior is empty~\citep{danao1994}. 
We expect the crossing of degenerate cones to induce non-smooth bifurcations because at the crossing of degenerate cones there is necessarily a continuum of solutions. Indeed, if $\det(C_{-M}(\alpha))=0$, the full orthant $\pos C_I(\alpha)$ is mapped by $f_{M}$ onto the (lower-dimensional) degenerate cone $\pos C_{-M}(\alpha)$. Thus, given $q\in \pos C_{-M}(\alpha)$, there must exist a (locally linear) subset of $\pos C_I(\alpha)$ that is mapped by $f_{M}$ to $q$~\citep{danao1994,murty1972}.

\begin{exmp}
  \label{exmp:lcp2:continuum}
  Let us consider the degenerate matrix 
  \begin{displaymath}
   M = \begin{bmatrix}
    1 & 1 \\ 1 & 1
  \end{bmatrix}
  \end{displaymath}
  and the path $\bar{q}_{a}$ as in \eqref{eq:lcp1:path:a}. For 
  $\alpha = \{1, 2\}$, solutions of \eqref{eq:lcp1:pwl}
  are characterized by the expression
  \begin{displaymath}
    \begin{bmatrix}
      4 \lambda - 4 \\ -4 \lambda
    \end{bmatrix} = 
    \begin{bmatrix}
      1 & 1 \\ 1 & 1
    \end{bmatrix} 
    \begin{bmatrix}
      x_{1} \\ x_{2}
    \end{bmatrix}, \quad x \in \pos C_{I}(\{1, 2\})
  \end{displaymath}
Note that the above equation has a nonempty solution set $S_{\{1, 2\}}$ if and only if 
$4 \lambda -4 = -4\lambda$, that is, if and only if $\lambda = \frac{1}{2}$. Hence, for 
$\lambda = \frac{1}{2}$ the solution set is given by
\begin{multline*}
  S_{\{1, 2\}} = \Big\{ (x, \lambda) \in \RE^{2} \times \Big\{\frac{1}{2}\Big\} 
    \mid \\ 
    x = \begin{bmatrix}
    \mu \\ -2 -\mu
\end{bmatrix}, \mu \in [-2, 0] \Big\},
\end{multline*}
whereas for the other subsets $\alpha \subset \{1, 2\}$, the solutions are
\begin{align*}
  S_{\emptyset} &= \emptyset
  \\
  S_{\{1\}} &= \left\{ (x, \lambda) \in \RE^{2} \times [0, 1/2) 
  \mid x = \begin{bmatrix} 4 \lambda - 4 \\ 4 - 8 \lambda
   \end{bmatrix} \right\}
   \\
  S_{\{2\}} &= \left\{ (x, \lambda) \in \RE^{2} \times (1/2, 1] \mid x = \begin{bmatrix}
  8 \lambda - 4 \\ -4 \lambda
\end{bmatrix}\right\}
\end{align*}

Therefore, for $\alpha = \{1, 2\}$ the solution set $S_{\alpha}$ has an infinite number of solutions
for a single value of $\lambda$, which corresponds to the situation in which the path $\bar{q}_{a}$ intersects
the degenerate cone $C_{M}(\alpha)$.
\end{exmp}

We summarize the results of this section as follows.
\begin{itemize}
 \item Non-smooth bifurcations can happen when the path defined by $q=\bar q(\lambda)$ crosses a face of non-degenerate cones, or at the crossing of degenerate cones.
 \item Crossing a degenerate cone always leads to bifurcations.
 \item The presence and nature of a bifurcation when crossing a face of a non-degenerate cone depends on the nature and disposition of the other cones that share that face.
\end{itemize}
 It follows that non-smooth bifurcations in LCPs are essentially determined by: i) the complementary cone configuration; ii) how the path moves across them. 

\section{Main results}
\label{section:main:results}

Similarly to smooth bifurcation theory, it is possible to use equivalence relations to provide an exhaustive list of the possible bifurcation phenomena. We start here this program by deriving a notion of \emph{equivalence between LCPs}, which will provide equivalence classes of cone configurations. The relevance of this notion in classifying non-smooth bifurcation problems will be then illustrated.

\subsection{Equivalence between cone configurations}

Our notion of equivalence between LCPs $(M,q)$ and $(N,r)$ has a topological and an algebraic 
component. The algebraic component captures the relations among the complementary cones
that $M$ and $N$ generate. The relevant algebraic structure is that of a Boolean algebra,
a subject that we now briefly recall (see~\cite{givant,sikorski} for more details).

Let $X$ be a set and $\mathcal{P}(X)$ the power set on $X$. A \emph{field of sets} is 
a pair $(X,\Fi)$ where $\Fi \subset \mathcal{P}(X)$ is closed under intersections of pairs 
of sets and complements of individual sets (this implies closure under union of pairs of sets).

Let $\Ge$ be a subset of $\mathcal{P}(X)$. The field of sets \emph{generated} by $\Ge$ is
the intersection of all the fields of sets that contain $\Ge$. 

A field of sets is a concrete example of a Boolean algebra, and as such, the usual algebraic
concepts apply to them.

\begin{defn}
 A \emph{Boolean homomorphism} from the field $(X,\Fi)$ onto the field $(X',\Fi')$ is a 
 mapping $h : \Fi \to \Fi'$ such that 
 \begin{equation} \label{eq:homo}
  h(P_1 \cap P_2) = h(P_1)\cap h(P_2) \quad \text{and} \quad h(-P_1) = -h(P_1)
 \end{equation}
 for all $P_1, P_2 \in \Fi$. Here, $-P_1$ denotes the complement of $P_1$. A one-to-one Boolean 
 homomorphism $h$ is called a \emph{Boolean isomorphism}. An isomorphism of a field 
 onto itself is called a \emph{Boolean automorphism}.
\end{defn}

\begin{defn}
 A Boolean mapping $h  : \Fi \to \Fi'$ is said to be \emph{induced} by a mapping 
 $\varphi : X' \to X$ if 
 \begin{equation} \label{eq:ind}
  h(P) = \varphi^{-1}(P)
 \end{equation}
 for every set $P \in \Fi$.
\end{defn}

Allow us to present a simple corollary to a theorem by Sikorski.

\begin{cor} \label{cor:ext}
 Let $\Fi$ be a field generated by $\Ge$. If a bijection $g : \Ge \to \Ge'$ is induced 
 by a bijection $\varphi : X' \to X$, then $g$ can be extended to a Boolean isomorphism
 $h : \Fi \to \Fi'$.
\end{cor}

\begin{pf}
 Define $h$ as in~\eqref{eq:ind}. Since $\varphi$ is bijection, $h$ satisfies~\eqref{eq:homo},
 that is, $g$ can been uniquely extended to a Boolean homomorphism from $\Fi$ into $\Fi'$. Likewise, $g^{-1}$
 can be extended to a Boolean homomorphism from $\Fi'$ into $\Fi$. It follows
 from~\cite[Thm. 12.1]{sikorski} that $h$ is indeed a Boolean isomorphism. \qed
\end{pf}

Now, consider the collection $\Ge_M = \left\{ \pos C_M(\alpha) \right\}_{\alpha}$, and 
let $(\RE^n,\Fi_M)$ be the field of sets generated by $\Ge_M$. We are now ready to state our main
definition.

\begin{defn}
 Two matrices $M, N \in \RE^{n\times n}$ are said to be \emph{LCP equivalent}, 
 $M \sim N$, if there exists topological isomorphisms (i.e., homeomorphisms)
 $\phi, \psi : \RE^n \to \RE^n$ such that 
 \begin{equation} \label{eq:topoEquiv}
  f_M = \varphi \circ f_N \circ \psi \;,
 \end{equation}
 where $\psi$ induces a Boolean automorphism on $\Fi_I$.
\end{defn}

Condition~\eqref{eq:topoEquiv} is the commutative diagram
\begin{displaymath}
\begin{tikzcd}
 \RE^n \arrow{r}{\psi} \arrow[swap]{d}{f_M} & \RE^n \arrow{d}{f_N} \\
 \RE^n                                      & \arrow{l}{\varphi} \RE^n
\end{tikzcd} \;.
\end{displaymath}
It is standard in the literature of singularity theory~\citep{arnold}, 
and ensures that we can continuously map solutions of the problem $f_M(x) = q$
into solutions of the problem $f_N(x') = \varphi^{-1}(q)$. The requirement on $\psi$ being a
Boolean automorphism implies that $\psi$ maps orthants into orthants, intersections of 
orthants into intersections of orthants, and so forth; and this ensures that the complementarity 
condition is not destroyed by the homeomorphisms.


\begin{thm} \label{thm:equiv}
 The matrices $M, N \in \RE^n$ are LCP equivalent if, and only if, there exists a bijection $g : \Ge_M \to \Ge_N$
 induced by a homeomorphism $\varphi : \RE^n \to \RE^n$.
\end{thm}

\begin{pf}
 Suppose that $M$ is LCP equivalent to $N$. Define $g$ as
 \begin{displaymath}
  g(\pos C_{M}(\alpha)) = \varphi^{-1}(\pos C_{M}(\alpha)) \;,
 \end{displaymath}
 where $\varphi$ satisfies~\eqref{eq:topoEquiv}. By~\eqref{eq:topoEquiv},
 \begin{displaymath}
  \varphi^{-1} \circ f_M(\pos C_I(\alpha)) = f_N \circ \psi(\pos C_I(\alpha)) \;.
 \end{displaymath}
 Since $\psi$ induces a Boolean automorphism on $\Fi_I$,
 \begin{displaymath}
  \psi(\pos C_I(\alpha)) = \pos C_I(\beta)
 \end{displaymath}
 for some $\beta$. Thus,
 \begin{displaymath}
  \varphi^{-1} \left( \pos C_{M}(\alpha) \right) = f_N \left( \pos C_I(\beta) \right) \;,
 \end{displaymath}
 so that
 \begin{equation} \label{eq:gBijec}
  g\left( \pos C_{M}(\alpha) \right) = \pos C_N(\beta) \;.
 \end{equation}
 This shows that $\varphi$ necessarily induces a bijection $g$ from $\Ge_M$ onto $\Ge_N$.

 For sufficiency, suppose that there is a bijection $g : \Ge_M \to \Ge_N$
 induced by a homeomorphism $\varphi : \RE^n \to \RE^n$. We will construct
 $\psi$ explicitly. Use the equation 
 \begin{displaymath}
  g(\pos C_{M}(\alpha)) = \pos C_{N}(\hat{\beta}(\alpha))
 \end{displaymath}
 to define the bijection $\hat{\beta}$ on the power set of $\{1,\dots,n\}$
 and denote its inverse by $\hat{\alpha}$. Now, define $\psi$ as
 \begin{multline*}
  \psi(x) = C_{-N}^{-1}(\hat{\beta}(\alpha))\cdot \varphi^{-1}\left( C_{-M}(\alpha) \cdot x \right) \;, \\
   \quad x \in \interior \pos C_I(\alpha) \;.
 \end{multline*}
 Note that $\psi(x) \in \pos C_I(\hat{\beta}(\alpha))$, so the application of $f_N$ on both sides of 
 the equation shows
 \begin{displaymath}
  f_N \circ \psi(x) = \varphi^{-1}\circ f_M(x)
 \end{displaymath}
 for $x$ in the interior of any orthant.

 Clearly, $\psi$ is piecewise continuous in the interior of the orthants. Its continuity at the boundaries
 follows from the continuity of $f_N$. More precisely,
 \begin{displaymath}
  x' = C_{-N}^{-1}(\beta_i) \cdot C_{-N}(\beta_j) \cdot x'
 \end{displaymath}
 for any indexes $\beta_i, \beta_j$ such that 
 \begin{displaymath}
  x' \in \pos C_{I}(\beta_i) \cap \pos C_{I}(\beta_j) \;.
 \end{displaymath}
 Thus, for $x$ in the boundary $\pos C_{I}(\alpha_i) \cap \pos C_{I}(\alpha_j)$, we have
 \begin{displaymath}
  C_{-N}^{-1}(\hat{\beta}(\alpha_i))\cdot \varphi^{-1}\left(f_M(x) \right) = 
  C_{-N}^{-1}(\hat{\beta}(\alpha_j))\cdot \varphi^{-1}\left(f_M(x) \right)
 \end{displaymath}
 so that the image of $x$ is the same, regardless of whether $\alpha_i$ or $\alpha_j$ is
 used in the definition of $\psi$.

 It is not difficult to verify that $\psi(x)$ is invertible with inverse
 \begin{multline*}
  \psi^{-1}(x') = C_{-M}^{-1}(\hat{\alpha}(\beta))\cdot \varphi\left( C_{-N}(\beta)\cdot x' \right) \;, \\
   x' \in \interior \pos C_I(\beta) \;.
 \end{multline*}
 Indeed, the composition $\psi^{-1}\circ \psi(x)$ gives
 \begin{multline*}
  C_{-M}^{-1}\left( \hat{\alpha}(\beta) \right)
   \varphi\left( C_{-N}(\beta) C_{-N}^{-1}(\hat{\beta}(\alpha)) \varphi^{-1}(C_{-M}(\alpha) x) \right) \;, \\
   x \in \interior \pos C_I(\alpha) \;.
 \end{multline*}
 This expression reduces to the identity by the fact that $\hat{\alpha}$ is the inverse of $\hat{\beta}$.
 
 By similar arguments, we can show that $\psi \circ \psi^{-1}$ is the identity, and that 
 $\psi^{-1}$ is continuous at the boundaries of the orthants. \qed
\end{pf}


\begin{rem} \label{rem:nec}
 It follows from Corollary~\ref{cor:ext} that a necessary condition for $M \sim N$ is the existence
 of a bijection $g : \Ge_M \to \Ge_N$ that extends to an isomorphism $h : \Fi_M \to \Fi_N$. 
\end{rem}

\begin{exmp} \label{ex:equiv}
Consider the matrices
\begin{displaymath}
 M = 
  \begin{bmatrix}
   -1 & 1 \\ 0.9 & -1
  \end{bmatrix}
  \;, \quad
 N = 
  \begin{bmatrix}
   -1 & 1 \\ 1.1 & -1
  \end{bmatrix}
  \quad \text{and} \quad
 O = 
  \begin{bmatrix}
   0.5 & 1 \\ 1 & 0.5
  \end{bmatrix}
  \;.
\end{displaymath}
Their cone configurations are shown in Fig.~\ref{fig:equivalenceEx}.

\begin{figure*}
\centering
 \includegraphics[width=0.7\textwidth]{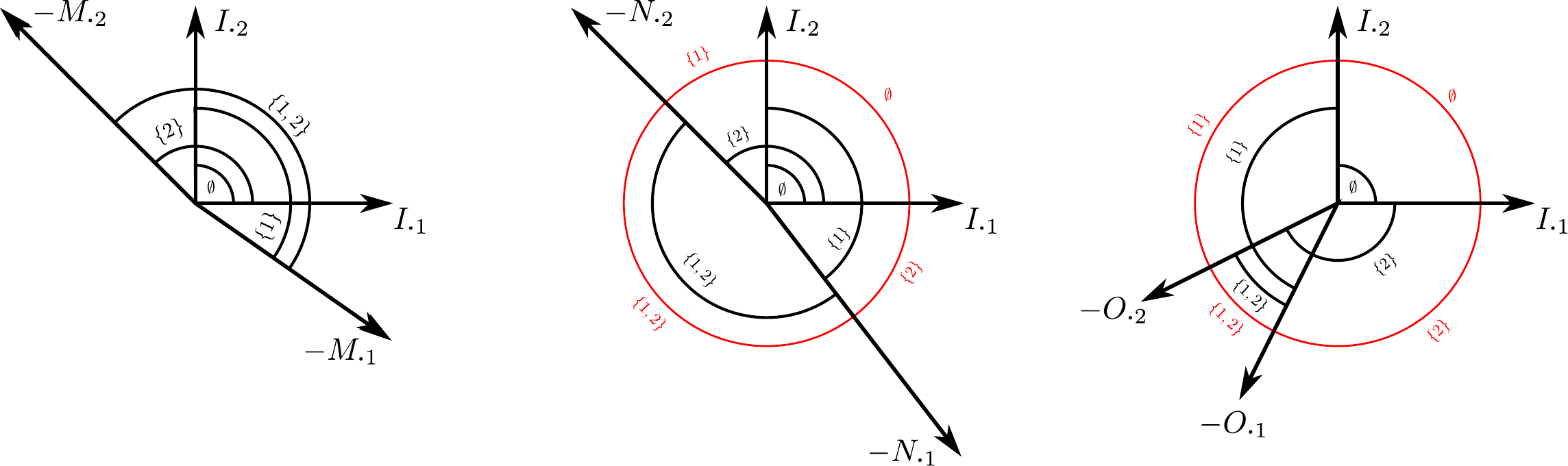}
\caption{Complementary cones of the matrices $M$, $N$ and $O$ in Example~\ref{ex:equiv}, depicted by
 black arcs. The cones generated by $C_{\bar{N}}(\alpha)$ and $C_{\bar{O}}(\alpha)$ are depicted by
 red arcs. The matrices $M$ and $N$ are not equivalent, but $N$ and $O$ are, as their complementary
 cones have the same Boolean structure.}
 \label{fig:equivalenceEx}
\end{figure*}

Note that 
\begin{displaymath}
 {\bigcap}_\alpha \pos C_M(\alpha) = \pos C_M(\emptyset) \;.
\end{displaymath}
Suppose, for the sake of argument, that there exists a bijection $g : \Ge_M \to \Ge_N$ 
that extends to an isomorphism $h : \Fi_M \to \Fi_N$. Then,
\begin{align*}
 h\left( {\bigcap}_\alpha \pos C_M(\alpha) \right) &= h\left( \pos C_M(\emptyset) \right) \\
    {\bigcap}_\alpha \pos C_N(\hat{\beta}(\alpha)) &= \pos C_N(\hat{\beta}(\emptyset))
\end{align*}
for some bijection $\hat{\beta}$. However, note that the intersection of all the complementary
cones generated by $N$ is no longer a cone. This is a contradiction, from which we conclude that such a $g$ cannot
exist and that, by Remark~\ref{rem:nec}, $M$ and $N$ are not equivalent. This is intuitively clear
since, depending on the location of $q$, there can be none, two, or four solutions to the
LCP $(M,q)$; whereas, depending on the location of $r$ there can be either one or three solutions to the LCP $(N,r)$.

Although $N$ and $O$ are fairly `distant' from each other, they are LCP equivalent. To see this, 
consider the matrices
\begin{displaymath}
 \bar{N} = 
  \begin{bmatrix}
   N_{\boldsymbol{\cdot} 2} & N_{\boldsymbol{\cdot} 1}
  \end{bmatrix}
  \quad \text{and} \quad
 \bar{O} = 
  \begin{bmatrix}
   O_{\boldsymbol{\cdot} 2} & O_{\boldsymbol{\cdot} 1}
  \end{bmatrix} \;.
\end{displaymath}
It is lengthy but straight forward to verify that the mapping $\varphi : \RE^2 \to \RE^2$, given by
\begin{displaymath}
 \varphi(y') = C_{\bar{N}}(\hat{\gamma}(\alpha))\cdot C_{\bar{O}}^{-1}(\alpha)\cdot y' \;, \quad y' \in \pos C_{\bar{O}}(\alpha)
\end{displaymath}
with $\hat{\gamma}(\emptyset) = \{1,2\}$, $\hat{\gamma}(\{1\}) = \{1\}$, $\hat{\gamma}(\{2\}) = \{ 2\}$ and
$\hat{\gamma}(\{1,2\}) = \emptyset$,
maps the cones $\pos C_{\bar{O}}(\alpha)$ to the cones $\pos C_{\bar{N}}(\hat{\gamma}(\alpha))$ 
(see Fig.~\ref{fig:equivalenceEx}). Clearly, $\varphi$ is a homeomorphism. Also, it induces the 
mapping $g : \Ge_N \to \Ge_O$ given by
\begin{displaymath}
 g\left( \pos C_N(\alpha) \right) = \pos C_O(-\alpha) \;.
\end{displaymath}
We have verified the conditions of Theorem~\ref{thm:equiv}. 
\end{exmp}

In the example, $M$ and $N$ are not equivalent, even though they are `close' to each other. 
This issue takes us to the following concept.

\begin{defn}
 A matrix $M \in \RE^{n\times n}$ is said to be \emph{LCP stable} if it is LCP equivalent to 
 every matrix that is sufficiently close to it.
\end{defn}




\subsection{Classification of LCPs on the plane}

The following results will provide a characterization of equivalence classes of stable matrices in $\RE^{2\times 2}$.

\begin{lem}\label{lem:planar:equivalence:suff}
	Let $M\in \RE^{2\times 2}$. If $M_{12},M_{21}\neq 0$ and $\det(M_{\alpha\alpha})\neq 0$, for all $\alpha\subseteq\{1,2\}$, then $M$ is stable.
\end{lem}

\begin{pf}
Let $E_M = \begin{bmatrix} I & -M \end{bmatrix}$, let $A_{\alpha}$ with $\alpha \subseteq \{1,2\}$
be the connected components of $\RE^2 - \bigcup_{i=1}^4 \pos E_{M\bscdot i}$, and note that
$A_\alpha$ are (not necessarily convex) cones that partition $\RE^2$, i.e.,
${\bigcup}_{\alpha}\bar{A}_\alpha = \RE^2$ and $A_\alpha \cap A_\beta = \emptyset$ for $\alpha \neq \beta$.
For every index $\alpha$, define $G_M(\alpha) \in \RE^{2\times2}$ as the submatrix of $E_M$ such that
\begin{displaymath}
 \boundary A_{\alpha} = \pos G_M(\alpha)_{\bscdot 1} \cup \pos G_M(\alpha)_{\bscdot 2}
\end{displaymath}
and $|G_M(\alpha)| > 0$\footnote{By the assumptions of the lemma, every $2\times2$-submatrix of $E_M$ is nonsingular,
so positivity of the determinant is ensured by suitably arranging the columns of $G_M(\alpha)$.}.
Let $M'=M+\varepsilon\tilde M $ be a perturbation of $M$. The matrix $E_{M'}$ varies smoothly as a function
of $\varepsilon$. In particular, $M'_{\bscdot i}\to M_{\bscdot i}$ as $\varepsilon\to 0$ (where the convergence is in some, 
and thus every, norm on $\RE^n$).

Let $\varphi : \RE^2 \to \RE^2$ be defined by
\[
 \varphi(y')= G_M(\alpha)\cdot G_{M'}^{-1}(\alpha)\cdot y' \;, \quad y' \in A'_\alpha
\]
with $G_{M'}(\alpha)$ the $\varepsilon$-perturbation of $G_M(\alpha)$ such that $|G_{M'}(\alpha)| > 0$.
Thus, $\varphi|_{A'_\alpha}$, $\alpha \subseteq \{1,2\}$, is continuous and bijective.  We now show that 
$\varphi$ is well-defined, continuous, and bijective at the cone boundaries, too. Let $A'_\alpha$ and $A'_\beta$ 
be two contiguous cones with common boundary $\pos E_{M' \bscdot i}$. Note that, if $y' \in A'_\alpha \cap A'_\beta$,
then $y' = \kappa E_{M' \bscdot i}$ for some $\kappa>0$, so that 
\begin{displaymath}
 \varphi|_{A'_\alpha}(y') =\kappa E_{M' \bscdot i} = \varphi|_{A_\beta}(y') \;.
\end{displaymath}
Thus $\varphi:\RE^2\to\RE^2$ is a homeomorphism. 
Moreover, since $\varphi(I_{\bscdot i})=I_{\bscdot i}$, $i=1,2$, and $\varphi(M'_{\bscdot i})=M_{\bscdot i}$, $i=1,2$, 
it follows that $\varphi^{-1}(\pos C_M(\alpha))=\pos C_{M'}(\alpha)$, that is, $\varphi^{-1}$ induces a bijection 
$g:\Ge_M \to \Ge_{M'}$. Then Theorem~\ref{thm:equiv} implies $M$ and $M'$ are equivalent.
\qed
\end{pf}

\begin{thm}\label{thm:planar equivalence}
	Two matrices $M, N \in \RE^{2\times 2}$ are equivalent if
	\begin{displaymath}
		M_{12} \cdot N_{12}>0,\ M_{21} \cdot N_{21} > 0
	\end{displaymath}
	and
	\begin{displaymath}
		\det(M_{\alpha \alpha})\cdot \det(N_{\alpha \alpha}) > 0 \;, \quad \alpha \subseteq \left\{ 1,2 \right\} \;.
	\end{displaymath}
\end{thm}

\begin{pf}
Under the hypothesis of the theorem, $M_t:=(1-t)M+tN$ is stable for all $t\in[0,1]$, because $M_t$ satisfies the conditions of Lemma \ref{lem:planar:equivalence:suff} for all $t\in[0,1]$. Thus, for all $t\in[0,1]$ there exists a neighborhood $U_t$ of $t$ in $[0,1]$ such that $M_t\sim M_{t'}$ for all $t'\in U_t$. Because $[0,1]$ is compact, we can cover it with a finite number of $U_t$'s, say, for $t_0<t_1<\cdots<t_n$, with $0\in U_{t_0}$ and $1\in U_{t_n}$. Let $\tau_i\in U_{t_{i-1}}\cap U_{t_i}$, $i=1,\ldots,n$. Then
$M=M_0\sim M_{\tau_1}\sim M_{\tau_2} \sim \cdots \sim M_{\tau_{n}}\sim M_1=N$.
\qed
\end{pf}

\begin{cor}\label{cor:planar equivalence neces}
	Let $M\in \RE^{2\times 2}$ with $\det(M_{\alpha\alpha})= 0$ for some $\alpha\subseteq\{1,2\}$, then $M$ is not stable.
\end{cor}

\begin{pf}
 Recalling that non-singular matrices are dense, there exists $M'$ arbitrarily close to $M$ such that $\det(M'_{\alpha\alpha})\neq 0$ for all $\alpha\subseteq\{1,2\}$. Then, any bijection $g:\mathcal G_M\to\mathcal G_{M'}$ cannot be induced by a homeomorphism $\varphi^{-1}$. If it was, $\varphi^{-1}$ would map an empty-interior complementarity cone to a non-empty interior complementarity cone, a contradiction. It follows by Theorem~\ref{thm:equiv} that $M$ cannot be stable. \qed
\end{pf}

The results of this section provides a list of ``normal forms" to explore equivalence classes of stable matrices in $\RE^{2\times 2}$. The matrices
\begin{displaymath}
	M_\delta = 
	\begin{bmatrix}
		\delta_1 & \delta_3 \\
		-\delta_3(2\delta_0-\delta_1\delta_2) & \delta_2 
	\end{bmatrix}
\end{displaymath}
and
\begin{displaymath}
N_\delta = 
\begin{bmatrix}
\delta_1 & \delta_3 \\
-\delta_3(0.5\delta_0-\delta_1\delta_2) & \delta_2  \;,
\end{bmatrix}
\end{displaymath}
with $\delta_i\in\{-1,1\}$, $i=0,\ldots,3$ span all possible 
sign combinations in the statement of Theorem~\ref{thm:planar equivalence}. Thus, any stable matrix satisfying the condition of the theorem is equivalent to $M_\delta$, for some combination of $\delta_i\in\{-1,1\}$, $i=0,\ldots,4$. By varying the parameters of $M_\delta$, we can thus construct an explicit list of equivalence classes of stable bidimensional matrices. The constructed list might not be exhaustive, but by Corollary~\ref{cor:planar equivalence neces} what is left out from this classification is the zero-measure set of matrices satisfying $\det(M_{\alpha\alpha})\neq 0$, for all $\alpha\subseteq\{1,2\}$, but $M_{12}M_{21}=0$. Stability and equivalence class of matrices in this zero-measure set are assessed {\it a posteriori} on a case-by-case basis in the normal form matrix
\begin{displaymath}
	O_\delta = 
	\begin{bmatrix}
		\delta_1 & \delta_3 \\
		\delta_4 & \delta_2 
	\end{bmatrix}
\end{displaymath}
with $\delta_1,\delta_2\in\{-1,1\}$ and $\delta_3,\delta_4\in\{-1,0,1\}$, $\delta_3\delta_4=0$.

\begin{figure}
\centering
\includegraphics[width=0.8\columnwidth]{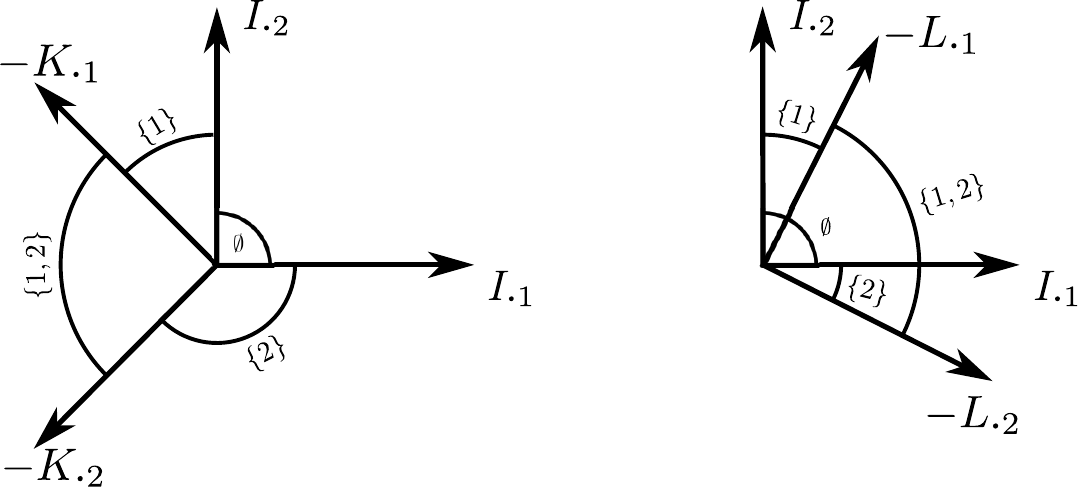}
\caption{Complementary cones of the matrices $K$ (left) and $L$ (right) defined in~\eqref{eq:KL},
  depicted by
 black arcs.}
 \label{fig:equivalenceEx2}
\end{figure}

After studying the cone structure of each of these matrices, we conclude that there are only four classes
of LCP stable matrices in $\RE^{2\times 2}$. Representative members of two different classes are the matrices
$M$ and $N$, defined in Example~\ref{ex:equiv}. Two more representative matrices are
\begin{equation} \label{eq:KL}
 K = 
	\begin{bmatrix}
	  1 & 1 \\
	 -1 & 1 
	\end{bmatrix}
 \quad \text{and} \quad
 L = 
	\begin{bmatrix}
	  -0.5 & -1 \\
	    -1 & 0.5 
	\end{bmatrix} \;.
\end{equation}
Since $\Ge_K$ partitions $\RE^2$ (see Fig.~\ref{fig:equivalenceEx2}), the LCP $(K,q)$ has a unique solution for every $q$. 
A matrix with this property is called a $P$-matrix~\citep{cottle}. The complementary cones of $L$ 
are also shown in Fig.~\ref{fig:equivalenceEx2}. Depending on $q$, the LCP $(L,q)$ may either have two or no solutions.


%
%

\subsection{Bifurcation realization via LCP interconnection}


The strong link between piecewise linear functions and LCP's, pointed out in Proposition \ref{prop:lcp:pwl:equivalence},
motivates us to restrict ourselves to piecewise linear paths through cone configurations. In this setting, the path itself
can be generated from the solution set of another LCP \citep{eaves1981, garcia1981}. This approach naturally leads us 
towards an interconnection framework reminiscent of circuit theory, 
in the sense that an intricate
high-dimensional LCP is treated as the result of the \emph{interconnection} of simpler LCP's.
Proceeding in this way we prove that, by selecting appropriate \emph{inputs} 
and \emph{outputs}, the \emph{feedback} interconnection of LCPs is again an LCP. 
Afterwards, we use this decomposition approach to obtain the unfoldings of the pitchfork singularity.

We start by considering two linear complementarity problems in their $z$-coordinates,
that is,
\begin{displaymath}
  \text{LCP} (M_{k}, \bar{q}_{k}) :
  \begin{cases}
    w_{k} = M_{k} z_{k} + \bar{q}_{k}
    \\
    \RE_{+}^{n_{k}} \ni w_{k} \perp z_{k} \in \RE_{+}^{n_{k}}
  \end{cases} \;,
\end{displaymath}
where $M_{k} \in \RE^{n_{k} \times n_{k}}$ and $\bar{q}_{k} \in \RE^{n_{k}}$, for 
$k \in \{a, b \}$. Let $z_{k} \in \RE^{n_{k}}$ be the \emph{output} of the $k$-th LCP 
and let $\bar{q}_{k} \in \RE^{n_{k}}$ take the role of \emph{input}. 
Additionally, consider the interconnection rule
\begin{equation}
  \label{eq:feedback:pattern}
  \bar{q}_{a} = H_{a} z_{b} + \bar{\theta}_{a} \;,
  \quad
  \bar{q}_{b} = H_{b} z_{a} + \bar{\theta}_{b} \;,
\end{equation}
where $H_{a} \in \RE^{n_{a} \times n_{b}}$, $H_{b} \in \RE^{n_{b} \times n_{a}}$ and 
$\bar{\theta}_{k} \in \RE^{n_{k}}$ are additional inputs available for further interconnection.
With this convention we have the following result.

\begin{prop}
  \label{prop:lcp:interconnection}
  The interconnection of linear complementarity problems under the pattern \eqref{eq:feedback:pattern} is again a linear complementarity problem. 
\end{prop}

\begin{pf}
  The interconnection of two LCPs under the pattern \eqref{eq:feedback:pattern} yields,
  \begin{equation}
    \label{eq:lcp:feedback}
    \begin{cases}
    \begin{bmatrix}
      w_{a} \\ w_{b}
    \end{bmatrix}
    = 
    \begin{bmatrix}
      M_{a} & H_{a}
      \\
      H_{b} & M_{b}
    \end{bmatrix}
    \begin{bmatrix}
      z_{a} \\ z_{b}
    \end{bmatrix}
    + 
    \begin{bmatrix}
      \bar{\theta}_{a} \\ \bar{\theta}_{b}
    \end{bmatrix}
    \\
    \RE_{+}^{n_{a} + n_{b}} \ni
    \begin{bmatrix}
      w_{a} \\ w_{b} 
    \end{bmatrix} \perp
    \begin{bmatrix}
      z_{a} \\ z_{b}
    \end{bmatrix} \in \RE_{+}^{n_{a} + n_{b}}
  \end{cases} \;.
  \end{equation}

The conclusion follows directly from \eqref{eq:lcp:feedback}, which is an LCP of 
dimension $n_{a} + n_{b}$ with extended input $\begin{bmatrix}
\bar{\theta}_{a} & \bar{\theta}_{b} \end{bmatrix}^{\top}$ and extended output
$\begin{bmatrix} z_{a} & z_{b} \end{bmatrix}^{\top}$.
\qed
\end{pf}
Note that, in contrast to the framework of dynamical systems, we are studying static 
relations that may be set-valued. Thus, the conditions for well-posedness of \eqref{eq:lcp:feedback} are more relaxed in comparison with their smooth counterpart.

\subsection{Realization of some non-smooth bifurcations and their unfolding: A non-smooth pleat}



\begin{figure}
  \centering
  \includegraphics[width=0.6\columnwidth]{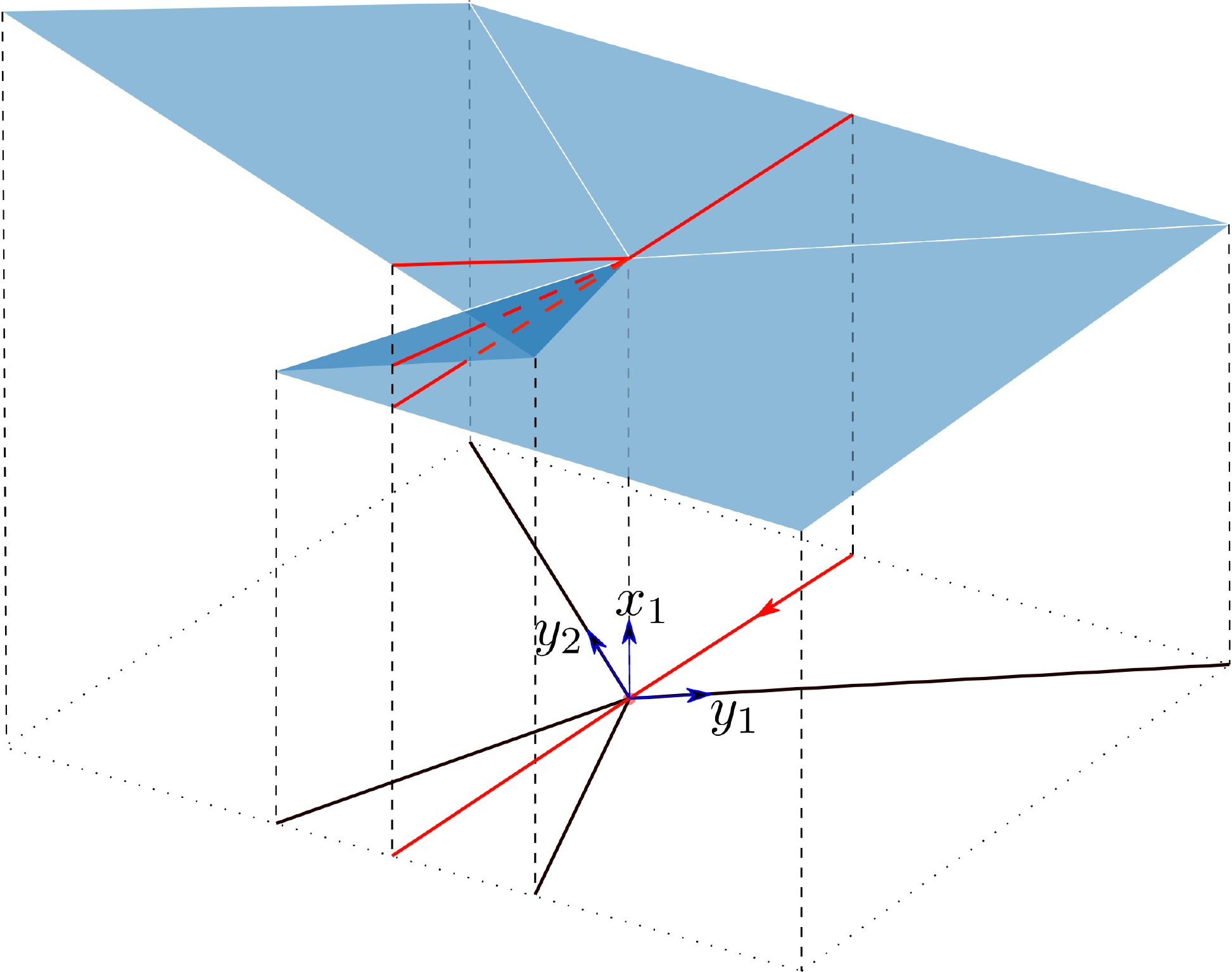}
  \caption{Non-smooth pleat, pitchfork path, and their projection to the plane 
  $(y_{1}, y_{2})$. The black lines on the plane are generators of complementary cones.}
  \label{fig:nonsmoothPleat}
\end{figure}

Let us consider the class of LCPs represented by the matrix $O$
in Fig.~\ref{fig:equivalenceEx}. This class gives rise to the non-smooth pleat
shown in Fig.~\ref{fig:nonsmoothPleat}. The pleat is given by
\begin{displaymath}
  \left\{ \begin{bmatrix} y_{1} & y_{2} & x_{1} \end{bmatrix}^{\top} \in \RE^{3} \mid
    \exists \; x_{2} \in \RE \text{ such that }  f_{O}(x) = y \right\} \;,
\end{displaymath}
where $f_{O}:\RE^{2} \to \RE^{2}$ is the piecewise linear map defined in Proposition 
\ref{prop:lcp:pwl:equivalence}.
It is worth to remark that the non-smooth pleat is stable in the sense that
the matrix $O$ is LCP-stable.

In complete analogy with the smooth case, see
e.g. \cite[Chapter III.12]{golubitsky1985}, one can recover a large family of bifurcations 
from the pleat by selecting appropriate paths through it. We illustrate this with the 
pitchfork singularity and its unfoldings, but it is also possible to obtain the 
hysteresis and the cusp singularities and their unfoldings by changing the path in a suitable way.

\begin{figure}
  \centering
  \includegraphics[width=0.8\columnwidth]{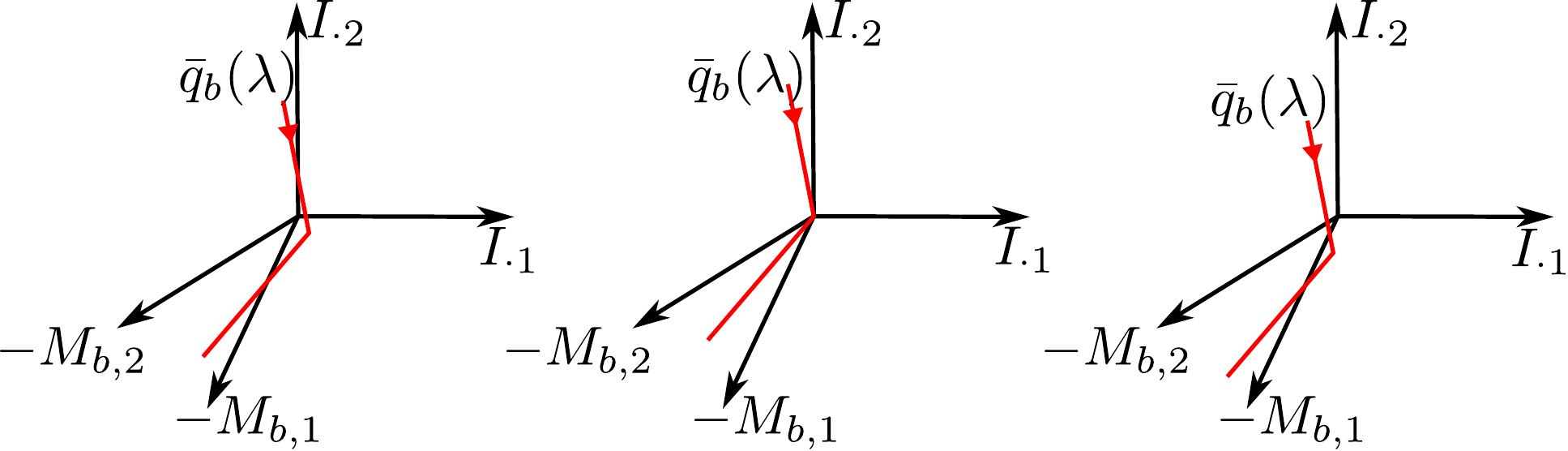}
  \caption{A sample of paths to recover the pitchfork singularity and its unfoldings.}
  \label{fig:pitchfork:paths}
\end{figure}

Concretely, let us consider the LCP $(M_{b}, \bar{q}_{b})$ associated to the non-smooth pleat shown in Fig.~\ref{fig:nonsmoothPleat}
with matrix $M_{b} = 2 O$ and $O$ as in Example 
\ref{ex:equiv}. In order to realize the path $\bar{q}_{b}$, we follow the interconnection approach described in the previous subsection.
We consider the second LCP $(M_{a}, \bar{q}_{a})$ with $M_{a} = 1$ and path $\bar{q}_{a}(\lambda) =
2 \lambda - 1$. The LCP $(M_{a}, \bar{q}_{a})$ has a unique solution for every $\lambda \in \RE$ which is computed easily as
\begin{equation}
  \label{eq:pleat:path}
  z_{a}(\lambda) = 
  \begin{cases}
    0, & \lambda < \frac{1}{2}
    \\
    2 \lambda - 1, & \frac{1}{2} \leq \lambda
  \end{cases} \;.
\end{equation}
We thus set the path $\bar{q}_{b}(\lambda)$ as
\begin{displaymath}
\bar{q}_{b}(\lambda) = R_{s} \begin{bmatrix}
    z_{a}(\lambda) \\ \lambda
  \end{bmatrix} + \begin{bmatrix}
    \mu_{1} \\ \mu_{2}
  \end{bmatrix} \;,
\end{displaymath}
where $R_{s}$ is a rotation matrix, $s$ is the angle of rotation and the parameters $\mu_{1}, \mu_{2}$ are extra degrees of freedom that
will allow us to change the path $\bar{q}_{b}(\lambda)$ on the pleat. Equivalently, the resulting LCP can be seen as the interconnection
between LCP $(M_{a}, \bar{q}_{a})$ and LCP $(M_{b}, \bar{q}_{b})$ under the interconnection rule \eqref{eq:feedback:pattern} with
\begin{align*}
  H_{a} & = 0 \;, & H_{b} & = \begin{bmatrix} \cos s \\ \sin s \end{bmatrix} \;, \\
  \bar{\theta}_{a} & = \mu_{1} - \lambda \sin s \;, & \bar{\theta}_{b} & = \mu_{2} + \lambda \cos s \;.
\end{align*}
Let us fix $s = \frac{10}{9} \pi$. By varying the parameters $\mu_{1}$ and $\mu_{2}$ we are able to displace
the path $\bar{q}_{b}(\lambda)$ on the pleat whose projection onto the plane is depicted in Fig.~\ref{fig:pitchfork:paths} for different
values of the vector $\mu$.

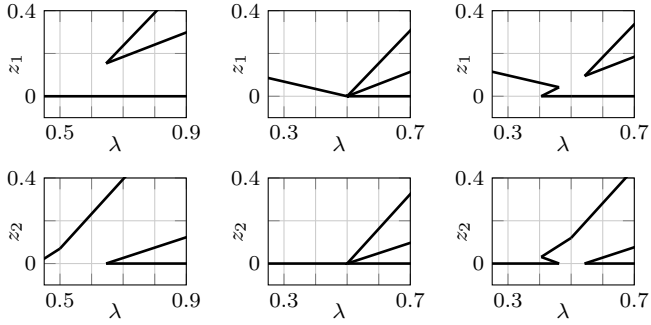
\begin{figure}
  \centering
  \small
  \begin{tikzpicture}
  \begin{groupplot}[group style={group size=1 by 2, vertical sep=8mm}, 
    height=3cm, width=0.39\columnwidth]
\nextgroupplot[
xtick = {0.5, 0.6, 0.7, 0.8, 0.9},
xticklabels = {$0.5$, , , , $0.9$},
ytick = {0, 0.2, 0.4},
yticklabels = {$0$, , $0.4$},
x grid style={white!80.0!black},
y grid style={white!80.0!black},
xmajorgrids,
ymajorgrids,
x label style={at={(axis description cs:0.5,-0.1)},anchor=north},
y label style={at={(axis description cs:-0.1,0.5)},anchor=south},
xlabel={$\lambda$},
ylabel={$z_{1}$},
xmin=0.45, xmax=0.9,
ymin=-0.1, ymax=0.4,
]
\addplot [line width=0.35mm, black, forget plot] table [x=l, y=x1, col sep=comma]{pitchfork_unfolding1_z_1.csv};
\addplot [line width=0.35mm, black, forget plot] table [x=l, y=x1, col sep=comma]{pitchfork_unfolding1_z_2.csv};
\addplot [line width=0.35mm, black, forget plot] table [x=l, y=x1, col sep=comma]{pitchfork_unfolding1_z_3.csv};
\addplot [line width=0.35mm, black, forget plot] table [x=l, y=x1, col sep=comma]{pitchfork_unfolding1_z_4.csv};
\addplot [line width=0.35mm, black, forget plot] table [x=l, y=x1, col sep=comma]{pitchfork_unfolding1_z_5.csv};
\addplot [line width=0.35mm, black, forget plot] table [x=l, y=x1, col sep=comma]{pitchfork_unfolding1_z_6.csv};

\nextgroupplot[
xtick = {0.5, 0.6, 0.7, 0.8, 0.9},
xticklabels = {$0.5$, , , , $0.9$},
ytick = {0, 0.2, 0.4},
yticklabels = {$0$, , $0.4$},
x grid style={white!80.0!black},
y grid style={white!80.0!black},
xmajorgrids,
ymajorgrids,
x label style={at={(axis description cs:0.5,-0.1)},anchor=north},
y label style={at={(axis description cs:-0.1,0.5)},anchor=south},
xlabel={$\lambda$},
ylabel={$z_{2}$},
xmin=0.45, xmax=0.9,
ymin=-0.1, ymax=0.4,
]
\addplot [line width=0.35mm, black, forget plot] table [x=l, y=x2, col sep=comma]{pitchfork_unfolding1_z_1.csv};
\addplot [line width=0.35mm, black, forget plot] table [x=l, y=x2, col sep=comma]{pitchfork_unfolding1_z_2.csv};
\addplot [line width=0.35mm, black, forget plot] table [x=l, y=x2, col sep=comma]{pitchfork_unfolding1_z_3.csv};
\addplot [line width=0.35mm, black, forget plot] table [x=l, y=x2, col sep=comma]{pitchfork_unfolding1_z_4.csv};
\addplot [line width=0.35mm, black, forget plot] table [x=l, y=x2, col sep=comma]{pitchfork_unfolding1_z_5.csv};
\addplot [line width=0.35mm, black, forget plot] table [x=l, y=x2, col sep=comma]{pitchfork_unfolding1_z_6.csv};

\end{groupplot}

\end{tikzpicture}
  \begin{tikzpicture}
  \begin{groupplot}[group style={group size=1 by 2, vertical sep=8mm}, 
    height=3cm, width=0.39\columnwidth]
\nextgroupplot[
xtick = {0.3, 0.4, 0.5, 0.6, 0.7},
xticklabels = {$0.3$, , , , $0.7$},
ytick = {0, 0.2, 0.4},
yticklabels = {$0$, , $0.4$},
x grid style={white!80.0!black},
y grid style={white!80.0!black},
xmajorgrids,
ymajorgrids,
x label style={at={(axis description cs:0.5,-0.1)},anchor=north},
y label style={at={(axis description cs:-0.1,0.5)},anchor=south},
xlabel={$\lambda$},
ylabel={$z_{1}$},
xmin=0.25, xmax=0.7,
ymin=-0.1, ymax=0.4,
]
\addplot [line width=0.35mm, black, forget plot] table [x=l, y=x1, col sep=comma]{pitchfork_orgCenter_z_1.csv};
\addplot [line width=0.35mm, black, forget plot] table [x=l, y=x1, col sep=comma]{pitchfork_orgCenter_z_2.csv};
\addplot [line width=0.35mm, black, forget plot] table [x=l, y=x1, col sep=comma]{pitchfork_orgCenter_z_3.csv};
\addplot [line width=0.35mm, black, forget plot] table [x=l, y=x1, col sep=comma]{pitchfork_orgCenter_z_4.csv};
\addplot [line width=0.35mm, black, forget plot] table [x=l, y=x1, col sep=comma]{pitchfork_orgCenter_z_5.csv};
\addplot [line width=0.35mm, black, forget plot] table [x=l, y=x1, col sep=comma]{pitchfork_orgCenter_z_6.csv};

\nextgroupplot[
xtick = {0.3, 0.4, 0.5, 0.6, 0.7},
xticklabels = {$0.3$, , , , $0.7$},
ytick = {0, 0.2, 0.4},
yticklabels = {$0$, , $0.4$},
x grid style={white!80.0!black},
y grid style={white!80.0!black},
xmajorgrids,
ymajorgrids,
x label style={at={(axis description cs:0.5,-0.1)},anchor=north},
y label style={at={(axis description cs:-0.1,0.5)},anchor=south},
xlabel={$\lambda$},
ylabel={$z_{2}$},
xmin=0.25, xmax=0.7,
ymin=-0.1, ymax=0.4,
]
\addplot [line width=0.35mm, black, forget plot] table [x=l, y=x2, col sep=comma]{pitchfork_orgCenter_z_1.csv};
\addplot [line width=0.35mm, black, forget plot] table [x=l, y=x2, col sep=comma]{pitchfork_orgCenter_z_2.csv};
\addplot [line width=0.35mm, black, forget plot] table [x=l, y=x2, col sep=comma]{pitchfork_orgCenter_z_3.csv};
\addplot [line width=0.35mm, black, forget plot] table [x=l, y=x2, col sep=comma]{pitchfork_orgCenter_z_4.csv};
\addplot [line width=0.35mm, black, forget plot] table [x=l, y=x2, col sep=comma]{pitchfork_orgCenter_z_5.csv};
\addplot [line width=0.35mm, black, forget plot] table [x=l, y=x2, col sep=comma]{pitchfork_orgCenter_z_6.csv};

\end{groupplot}

\end{tikzpicture}
  \begin{tikzpicture}
  \begin{groupplot}[group style={group size=1 by 2, vertical sep=8mm}, 
    height=3cm, width=0.39\columnwidth]
\nextgroupplot[
xtick = {0.3, 0.4, 0.5, 0.6, 0.7},
xticklabels = {$0.3$, , , , $0.7$},
ytick = {0, 0.2, 0.4},
yticklabels = {$0$, , $0.4$},
x grid style={white!80.0!black},
y grid style={white!80.0!black},
xmajorgrids,
ymajorgrids,
x label style={at={(axis description cs:0.5,-0.1)},anchor=north},
y label style={at={(axis description cs:-0.1,0.5)},anchor=south},
xlabel={$\lambda$},
ylabel={$z_{1}$},
xmin=0.25, xmax=0.7,
ymin=-0.1, ymax=0.4,
]
\addplot [line width=0.35mm, black, forget plot] table [x=l, y=x1, col sep=comma]{pitchfork_unfolding2_z_1.csv};
\addplot [line width=0.35mm, black, forget plot] table [x=l, y=x1, col sep=comma]{pitchfork_unfolding2_z_2.csv};
\addplot [line width=0.35mm, black, forget plot] table [x=l, y=x1, col sep=comma]{pitchfork_unfolding2_z_3.csv};
\addplot [line width=0.35mm, black, forget plot] table [x=l, y=x1, col sep=comma]{pitchfork_unfolding2_z_4.csv};
\addplot [line width=0.35mm, black, forget plot] table [x=l, y=x1, col sep=comma]{pitchfork_unfolding2_z_5.csv};
\addplot [line width=0.35mm, black, forget plot] table [x=l, y=x1, col sep=comma]{pitchfork_unfolding2_z_6.csv};

\nextgroupplot[
xtick = {0.3, 0.4, 0.5, 0.6, 0.7},
xticklabels = {$0.3$, , , , $0.7$},
ytick = {0, 0.2, 0.4},
yticklabels = {$0$, , $0.4$},
x grid style={white!80.0!black},
y grid style={white!80.0!black},
xmajorgrids,
ymajorgrids,
x label style={at={(axis description cs:0.5,-0.1)},anchor=north},
y label style={at={(axis description cs:-0.1,0.5)},anchor=south},
xlabel={$\lambda$},
ylabel={$z_{2}$},
xmin=0.25, xmax=0.7,
ymin=-0.1, ymax=0.4,
]
\addplot [line width=0.35mm, black, forget plot] table [x=l, y=x2, col sep=comma]{pitchfork_unfolding2_z_1.csv};
\addplot [line width=0.35mm, black, forget plot] table [x=l, y=x2, col sep=comma]{pitchfork_unfolding2_z_2.csv};
\addplot [line width=0.35mm, black, forget plot] table [x=l, y=x2, col sep=comma]{pitchfork_unfolding2_z_3.csv};
\addplot [line width=0.35mm, black, forget plot] table [x=l, y=x2, col sep=comma]{pitchfork_unfolding2_z_4.csv};
\addplot [line width=0.35mm, black, forget plot] table [x=l, y=x2, col sep=comma]{pitchfork_unfolding2_z_5.csv};
\addplot [line width=0.35mm, black, forget plot] table [x=l, y=x2, col sep=comma]{pitchfork_unfolding2_z_6.csv};

\end{groupplot}

\end{tikzpicture}

  \caption{The pitchfork singularity and its unfoldings, obtained from the paths depicted in Fig.~\ref{fig:pitchfork:paths}.}
  \label{fig:pitchfork}
\end{figure}

The associated bifurcation diagrams to the paths on Fig.~\ref{fig:pitchfork:paths} are shown in Fig.~\ref{fig:pitchfork}.
Note that the central path in Fig.~\ref{fig:pitchfork:paths} produces the pitchfork organizing center, whereas perturbations of this 
path lead to any of the left or right-hand side diagrams.

\section{Discussion and future directions}
\label{section:discussion:future}

We have presented a notion of global equivalence between LCPs that allows us to make a classification of this problems in the planar case.
In addition, an interconnection approach for the realization of non-smooth bifurcations was presented. These tools are thought to be handful for many applications,
as for instance, the analysis and design of neuromorphic circuits \citep{castanos2017}, the study of economic equilibria in competitive markets \citep{nagurney1993},
and the analysis of elastic-plastic structures in engineering \citep{pang1979}, just to name a few. 
This work also opens the path towards the analysis of behaviors in dynamical linear complementarity systems \citep{schaft1998}.





\bibliography{lcp}             

\begin{thebibliography}{22}
\providecommand{\natexlab}[1]{#1}
\providecommand{\url}[1]{\texttt{#1}}
\providecommand{\urlprefix}{URL }
\expandafter\ifx\csname urlstyle\endcsname\relax
  \providecommand{\doi}[1]{doi:\discretionary{}{}{}#1}\else
  \providecommand{\doi}{doi:\discretionary{}{}{}\begingroup
  \urlstyle{rm}\Url}\fi

\bibitem[{Acary et~al.(2011)Acary, Bonnefon, and Brogliato}]{acary2011}
Acary, V., Bonnefon, O., and Brogliato, B. (2011).
\newblock \emph{Nonsmooth modeling and simulation for switched circuits}.
\newblock Lecture Notes in Electrical Engineering. Springer.

\bibitem[{Arnold et~al.(1985)Arnold, Varchenko, and Gusein-Zade}]{arnold}
Arnold, V.I., Varchenko, A.N., and Gusein-Zade, S.M. (1985).
\newblock \emph{Singularities of Differentiable Maps}, volume~1.
\newblock Birkh{\"a}user.

\bibitem[{Brogliato(1999)}]{brogliato1999}
Brogliato, B. (1999).
\newblock \emph{Nonsmooth mechanics: models, dynamics and control}.
\newblock Springer-Verlag, London, second edition.

\bibitem[{Casta\~{n}os and Franci(2017)}]{castanos2017}
Casta\~{n}os, F. and Franci, A. (2017).
\newblock Implementing robust neuromodulation in neuromorphic circuits.
\newblock \emph{Neurocomputing}, 233, 3--13.

\bibitem[{Clarke(1990)}]{clarke1990}
Clarke, F. (1990).
\newblock \emph{Optimization and nonsmooth analysis}.
\newblock Society for Industrial and Applied Mathematics.

\bibitem[{Cottle et~al.(2009)Cottle, Pang, and Stone}]{cottle}
Cottle, R.W., Pang, J.S., and Stone, R.E. (2009).
\newblock \emph{The Linear Complementarity Problem}.
\newblock Classics in Applied Mathematics. Society for Industrial and Applied
  Mathematics.

\bibitem[{Danao(1994)}]{danao1994}
Danao, R.A. (1994).
\newblock Q-matrices and boundedness of solutions to linear complementarity
  problems.
\newblock \emph{Journal of Optimization Theory and Applications}, 83(2),
  321--332.

\bibitem[{Di~Bernardo et~al.(2008)Di~Bernardo, Budd, Champneys, Kowalczyk,
  Nordmark, Tost, and Piiroinen}]{diBernardo2008}
Di~Bernardo, M., Budd, C.J., Champneys, A.R., Kowalczyk, P., Nordmark, A.B.,
  Tost, G.O., and Piiroinen, P.T. (2008).
\newblock Bifurcations in nonsmooth dynamical systems.
\newblock \emph{SIAM review}, 50(4), 629--701.

\bibitem[{Dontchev and Rockafellar(2014)}]{dontchev2014}
Dontchev, A.L. and Rockafellar, R.T. (2014).
\newblock \emph{Implicit Functions and Solution Mappings}.
\newblock Springer, New York, second edition.

\bibitem[{Eaves and Lemke(1981)}]{eaves1981}
Eaves, B.C. and Lemke, C.E. (1981).
\newblock Equivalence of {LCP} and {PLS}.
\newblock \emph{Mathematics of Operation Research}, 6(4), 475--484.

\bibitem[{Garcia et~al.(1983)Garcia, Gould, and Turnbull}]{garcia1981}
Garcia, C.B., Gould, F.J., and Turnbull, T.R. (1983).
\newblock Relations between {PL} maps, complementary cones, and degree in
  linear complementarity problems.
\newblock In \emph{Homotopy Methods and Global Convergence}, 91--144. Plenum
  Publishing Corporation.

\bibitem[{Givant and Halmos(2009)}]{givant}
Givant, S. and Halmos, P.R. (2009).
\newblock \emph{Introduction to Boolean Algebras}.
\newblock Springer-Verlag, New York.

\bibitem[{Golubitsky and Schaeffer(1985)}]{golubitsky1985}
Golubitsky, M. and Schaeffer, D. (1985).
\newblock \emph{Singularities and Groups in Bifurcation Theory}, volume~I of
  \emph{Applied Mathematical Sciences}.
\newblock Springer, New York.

\bibitem[{Leenaerts and Bokhoven(1998)}]{leenaerts1998}
Leenaerts, D.M. and Bokhoven, W.M.G. (1998).
\newblock \emph{Piecewise linear modeling and analysis}.
\newblock Kluwer Academic. Springer, New York.

\bibitem[{Leine and Nijmeijer(2004)}]{leine2004}
Leine, R.I. and Nijmeijer, H. (2004).
\newblock \emph{Dynamics and bifurcations on non-smooth mechanical systems}.
\newblock Springer, Berlin.

\bibitem[{Murty(1988)}]{murty1988}
Murty, K.G. (1988).
\newblock \emph{Linear complementarity, linear and nonlinear programming}.
\newblock Helderman Verlag, Berlin.

\bibitem[{Murty(1972)}]{murty1972}
Murty, K.G. (1972).
\newblock On the number of solutions to the complementarity problem and
  spanning properties of complementary cones.
\newblock \emph{Linear Algebra and Its Applications}, 5(1), 65--108.

\bibitem[{Nagurney(1999)}]{nagurney1993}
Nagurney, A. (1999).
\newblock \emph{Network economics: A variational inequality approach}.
\newblock Advances in Computational Economics. Springer-Science+Business Media.

\bibitem[{Pang et~al.(1979)Pang, Kaneko, and Hallman}]{pang1979}
Pang, J., Kaneko, I., and Hallman, W. (1979).
\newblock On the solution of some (parametric) linear complementarity problems
  with applications to portfolio selection, structural engineering and
  actuarial graduation.
\newblock \emph{Mathematical Programming}, 16, 325 -- 347.

\bibitem[{Sikorski(1969)}]{sikorski}
Sikorski, R. (1969).
\newblock \emph{Boolean Algebras}.
\newblock Springer-Verlag, Berlin, Germany.

\bibitem[{Simpson(2010)}]{simpson2010}
Simpson, D.J.W. (2010).
\newblock \emph{Bifurcations in piecewise-smooth continuous systems}.
\newblock World Scientific, Singapore.

\bibitem[{van~der Schaft and Schumacher(1998)}]{schaft1998}
van~der Schaft, A.J. and Schumacher, J. (1998).
\newblock Complementarity modeling of hybrid systems.
\newblock \emph{{IEEE} Trans. Autom. Control}, 43, 483 -- 490.

\end{thebibliography}

\end{document}